\documentclass[10pt]{article}
\usepackage{amssymb, amsmath}
\usepackage{float,epsfig}

\begin{document}

\newtheorem{theorem}{\bf Theorem}[section]
\newtheorem{proposition}[theorem]{\bf Proposition}
\newtheorem{definition}[theorem]{\bf Definition}
\newtheorem{corollary}[theorem]{\bf Corollary}
\newtheorem{example}[theorem]{\bf Example}
\newtheorem{exam}[theorem]{\bf Example}
\newtheorem{remark}[theorem]{\bf Remark}
\newtheorem{lemma}[theorem]{\bf Lemma}
\newcommand{\nrm}[1]{|\!|\!| {#1} |\!|\!|}

\newcommand{\ba}{\begin{array}}
\newcommand{\ea}{\end{array}}
\newcommand{\von}{\vskip 1ex}
\newcommand{\vone}{\vskip 2ex}
\newcommand{\vtwo}{\vskip 4ex}
\newcommand{\dm}[1]{ {\displaystyle{#1} } }

\newcommand{\be}{\begin{equation}}
\newcommand{\ee}{\end{equation}}
\newcommand{\beano}{\begin{eqnarray*}}
\newcommand{\eeano}{\end{eqnarray*}}
\newcommand{\inp}[2]{\langle {#1} ,\,{#2} \rangle}
\def\bmatrix#1{\left[ \begin{matrix} #1 \end{matrix} \right]}
\def \noin{\noindent}
\newcommand{\evenindex}{\Pi_e}



\def \R{{\mathbb R}}
\def \C{{\mathbb C}}
\def \K{{\mathbb K}}
\def \J{{\mathbb J}}
\def \Lb{\mathrm{L}}

\def \T{{\mathbb T}}
\def \Pb{\mathrm{P}}
\def \N{{\mathbb N}}
\def \Ib{\mathrm{I}}
\def \Ls{{\Lambda}_{m-1}}
\def \Gb{\mathrm{G}}
\def \Hb{\mathrm{H}}
\def \Lam{{\Lambda_{m}}}
\def \Qb{\mathrm{Q}}
\def \Rb{\mathrm{R}}
\def \Mb{\mathrm{M}}
\def \norm{\nrm{\cdot}\equiv \nrm{\cdot}}

\def \P{{\mathbb{P}}_m(\C^{n\times n})}
\def \A{{{\mathbb P}_1(\C^{n\times n})}}
\def \H{{\mathbb H}}
\def \L{{\mathbb L}}
\def \G{{\F_{\tt{H}}}}
\def \S{\mathbb{S}}
\def \sigmin{\sigma_{\min}}
\def \elam{\Lambda_{\epsilon}}
\def \slam{\Lambda^{\S}_{\epsilon}}
\def \Ib{\mathrm{I}}
\def \Tb{\mathrm{T}}
\def \d{{\delta}}

\def \Lb{\mathrm{L}}
\def \N{{\mathbb N}}
\def \Ls{{\Lambda}_{m-1}}
\def \Gb{\mathrm{G}}
\def \Hb{\mathrm{H}}
\def \Delta{\triangle}
\def \Rar{\Rightarrow}
\def \p{{\mathsf{p}(\lam; v)}}

\def \D{{\mathbb D}}

\def \tr{\mathrm{Tr}}
\def \cond{\mathrm{cond}}
\def \lam{\lambda}
\def \sig{\sigma}
\def \sign{\mathrm{sign}}

\def \ep{\epsilon}
\def \diag{\mathrm{diag}}
\def \rev{\mathrm{rev}}
\def \vec{\mathrm{vec}}

\def \sk{\mathsf{skew}}
\def \sy{\mathsf{sym}}
\def \en{\mathrm{even}}
\def \odd{\mathrm{odd}}
\def \rank{\mathrm{rank}}
\def \pf{{\bf Proof: }}
\def \dist{\mathrm{dist}}
\def \rar{\rightarrow}

\def \rank{\mathrm{rank}}
\def \pf{{\bf Proof: }}
\def \dist{\mathrm{dist}}
\def \Re{\mathsf{Re}}
\def \Im{\mathsf{Im}}
\def \re{\mathsf{re}}
\def \im{\mathsf{im}}

\def \sym{\mathsf{sym}}
\def \sksym{\mathsf{skew\mbox{-}sym}}
\def \odd{\mathrm{odd}}
\def \even{\mathrm{even}}
\def \herm{\mathsf{Herm}}
\def \skherm{\mathsf{skew\mbox{-}Herm}}
\def \str{\mathrm{ Struct}}
\def \eproof{$\blacksquare$}
\def \proof{\noin\pf}

\def \bS{{\bf S}}
\def \cA{{\cal A}}
\def \E{{\mathcal E}}
\def \X{{\mathcal X}}
\def \F{{\mathcal F}}
\def \cH{\mathcal{H}}
\def \cJ{\mathcal{J}}
\def \tr{\mathrm{Tr}}
\def \range{\mathrm{Range}}
\def \adj{\star}

\def \pal{\mathrm{palindromic}}
\def \palpen{\mathrm{palindromic~~ pencil}}
\def \palpoly{\mathrm{palindromic~~ polynomial}}
\def \odd{\mathrm{odd}}
\def \even{\mathrm{even}}


\title{\sc  Structured mapping problems for linearly structured matrices}
\author{ Bibhas Adhikari\thanks{CoE Systems Science,
IIT Jodhpur, India, E-mail:
bibhas@iitj.ac.in } \, and Rafikul Alam\thanks{Department of Mathematics,
IIT Guwahati, India, E-mail: rafik@iitg.ernet.in, Fax:~+91-361-2690762/2582649. }  }
\date{}

\maketitle
\thispagestyle{empty}

{\small \noin{\bf Abstract.}  Given an appropriate class of  structured matrices $\S,$   we characterize matrices $X$ and $B$  for which  there exists a matrix $A \in \S$ such that $AX = B$ and determine all matrices in $\S$ mapping $X$ to $B.$ We also determine all matrices in $\S$ mapping $X$ to $B$ and having the smallest norm. We use these results to investigate structured backward errors of approximate eigenpairs and approximate invariant subspaces, and structured pseudospectra of structured matrices.}

\vone \noin{\bf Keywords.} Structured matrices, structured backward errors, Jordan and Lie algebras,   eigenvalues, eigenvectors, invariant subspaces.

\vone\noin {\bf AMS subject classification(2000):}  15A24,  65F15, 65F18,  15A18.

\section{Introduction} Consider a stable linear time-invariant (LTI) control  system
\be \label{paslti}
\begin{array}{lcl}
  \dot x & = & A x +B u,\ x(0)=0, \\
  y & = & C x +D u,
\end{array}
\ee
with   $A\in \mathbb K^{n\times n}$,
$B\in \mathbb K^{n\times p}$, $C\in \mathbb K^{p\times n}$ and
$D\in \mathbb K^{p\times p}.$ Here $ \K := \R \mbox{ or } \C,$  $u$ is the input, $x$ is the state and $y$ is the output.
The system (\ref{paslti}) is said to be passive if the  {\em Hamiltonian matrix}
\begin{equation}\label{hampas}
\cH=\bmatrix{ F & G \\ H & -F^*}  :=\bmatrix{ A-BR^{-1} C & -BR^{-1}B^* \\
-C^*  R^{-1} C &  -(A-BR^{-1}C)^*}
\end{equation}
has no purely imaginary eigenvalues, where  $R := D +
D^* $, see~\cite{anto,griv, AlaBKMM}.  A matrix $\cH\in {\mathbb K}^{2n\times 2n}$
of the form $\cH=\bmatrix{ A & F \\ G & -A^*}$
is called  Hamiltonian, where $G^* =G$ and $ F^*=F.$
Equivalently, $\cH$ is Hamiltonian $\iff  (\cJ\cH)^* = \cJ \cH,$ where
    $\cJ := \bmatrix{ 0 & I \\ -I & 0}$ and   $I$  the  identity matrix of size $n.$

For passivation problem, when purely imaginary eigenvalues occur, one tries to perturb $\cH$  by a Hamiltonian matrix $\E$ with {\em small norm} so that the perturbed matrix  $\cH+ \E$  has no purely imaginary eigenvalues. If such an $\E$ exists, then for some $ X \in \K^{2n\times p}$ and $ D \in  \K^{p\times p},$ we have
 $$ (\cH+\E)X = X D \Longrightarrow \E X = B :=\cH X-XD.$$ This leads us to the following mapping problem.

 \vone
 {\sc Problem~1.} {\bf (Hamiltonian mapping problem)} Given  $ X , B \in \K^{2n \times p},$ consider
 \beano  \mathrm{Ham}(X,B) &:=& \{ \cH \in \K^{2n\times 2n} : (\cJ \cH)^* = \cJ \cH \mbox{ and } \cH X = B\}, \\
  \sigma^{\rm Ham}(X, B) &:= & \inf\{ \|\cH\| : \cH \in  \mathrm{Ham}(X, B) \}. \eeano
  \begin{itemize} \item Characterize $X, B \in \K^{2n\times p} $ for which $\mathrm{Ham}(X, B) \neq \emptyset$  and determine all matrices in $\mathrm{Ham}(X, B).$
\item Also determine all optimal solutions $\cH_o \in \mathrm{Ham}(X, B)$ such that $\|\cH_o  \| = \sigma^{\rm Ham}(X, B).$
\end{itemize}

Motivated by {\sc Problem~1}, we now consider structured mapping problem for various classes of structured matrices. Let $\S$ denote a class of structured matrices in $\K^{n\times n}.$ The class $\S$ we consider in this paper is either a Jordan or a Lie algebra associated with an appropriate scalar product on $\K^n.$ This provides a general setting that encompasses important classes of structured matrices such as Hamiltonian, skew-Hamiltonian, symmetric, skew-symmetric, pseudosymmetric, persymmetric, Hermitian, Skew-Hermitian, pseudo-Hermitian, pseudo-skew-Hermitian, to name only a few, see~\cite{fran6}. We, therefore, consider the following problem.

\vone
\noin{\sc Problem~2.} {\bf (Structured Mapping Problem)} Let $\S \subset \K^{n\times n}$ be a class of structured matrices and let $ X, B \in \K^{n\times p}.$
Set \beano  \S(X,B) &:=& \{ A \in \S : A X = B\}, \\
  \sigma^{\S}(X, B) &:= & \inf\{ \|A\| : A \in  \S(X, B) \}. \eeano
 \begin{itemize} \item {\bf Existence:} Characterize $X, B \in \K^{n\times p} $ for which $\S(X, B) \neq \emptyset.$
 \item {\bf Characterization:} Determine all matrices in $\S(X, B).$ Also determine all optimal solutions $A_o \in \S(X, B)$ such that $\|A_o\| = \sigma^{\S}(X, B).$
\end{itemize}

We  mention that structured backward error of an approximate invariant subspace of a structured matrix also leads to a structured mapping problem.
A subspace   $ \mathcal{ X}$  is invariant under $A$ if $ A \mathcal{X} \subset \mathcal{X}.$

\vone
\noin
{\sc  Problem~3.} {\bf (Structured backward error)} Let  $ \S \subset \K^{n\times n}$ be a class of structured matrices and $A \in \S.$ Let $\mathcal{X}$ be a subspace of $\K^n.$  Set
$$\omega^{\S}(A, \mathcal{X}) := \min\{ \|\Delta A\| : \Delta A \in \S \mbox{ and } (A+\Delta A)\mathcal{X} \subset \mathcal{X}\}.$$
Find all $ E \in \S$ such $(A+E) \mathcal{X} \subset \mathcal{X}$ and $\|E\| = \omega^\S(A, \mathcal{X}).$

\vone

 If such a matrix  $E \in \S$ exists then  $(A+E) \mathcal{X} \subset \mathcal{X} \Rightarrow  EV = VR - AV =: B$ for some $R$ and a full column rank matrix $V$ whose columns form a basis of ${\cal X}.$ This shows that structured mapping problem naturally arises when analyzing structured backward error of an approximate invariant subspace.

Solutions of structured and unstructured mapping problems for a pair of vectors $x$ and $b$ in $\K^n$ have been studied extensively, see~\cite{fran6} and the references therein. In fact, for a pair of vectors $x$ and $b$ in $\K^n,$ a complete solution of the structured mapping problem has been provided in~\cite{fran6} when the class $\S \subset \K^{n\times n}$ of structured matrices is a Jordan or a Lie algebra associated with an orthosymmetric scalar product on $\K^n.$  For a pair of matrices $X$ and $B $ in $\K^{n\times p},$ existence and characterization of solutions to  $AX =B$ have been  discussed for  Hermitian solutions in~\cite{khatri, sun93} and for skew-Hermitian and symmetric solutions in~\cite{tiss03}. Also, for the Frobenius norm, \cite{sun93} provides an optimal Hermitian solution and \cite{tiss03} provides optimal skew-Hermitian and symmetric solutions to $AX=B.$

The main contributions of this paper are as follows. We provide a complete solution of the structured mapping problem ({\sc Problem~2}) when the class $\S \subset \K^{n\times n}$ of structured matrices is a Jordan or a Lie algebra associated with an orthosymmetric scalar product on $\K^n.$  We show that for the spectral norm there are infinitely many optimal solutions whereas for the Frobenius norm the optimal solution is unique. We determine all optimal solutions for the spectral and the Frobenius norms. We show that the results in~\cite{fran6} obtained for a pair of vectors follow as special cases of our general results. Finally, as an application of the structured mapping problem, we analyze structured backward errors of approximate invariant subspaces,  approximate eigenpairs, and structured pseudospectra of structured matrices.

\vone

\noindent{\bf Notation.} We denote the spectrum and the trace of a square matrix $A$ by $\Lambda(A)$ and $\tr(A),$ respectively. Also, we denote the subspace spanned by the columns of an $m$-by-$n$ matrix $X$ by $\mathrm{span}(X)$ and the complex conjugate of $X$ by $\bar{X}.$

\section{Linearly structured matrices} Let $\K^{n\times n}$ denote the set of all $n$-by-$n$ matrices with entries in $\K,$ where $\K = \R$ or $\C.$
We denote the transpose of a matrix $X \in \K^{m\times n}$ by $X^T$ and the conjugate transpose by $X^H.$

We now briefly define structured matrices that we consider in this paper; see~\cite{fran6} for further details.
Let  $ M \in \K^{n\times n}$ be unitary. Assume further that $M$ is either symmetric  or skew-symmetric or Hermitian or skew-Hermitian.
Define the scalar product $\inp{\cdot}{\cdot}_M : \K^n\times \K^n \rar \K$ by \be \label{sprod} \inp{x}{y}_M := \left\{\ba{cc} y^TMx, & \mbox{ bilinear form, }\\
y^HMx, & \mbox{ sesquilinear form.} \ea\right.\ee Then for $A \in \K^{n\times n}$ there is a unique adjoint operator $A^\star$ relative to the scalar product (\ref{sprod}) such that
 $\inp{Ax}{y}_M= \inp{x}{A^\star y}_M $ for all $x, y \in \K^n.$ The adjoint $A^\star$ is explicitly given by
\be \label{adj} A^{\star} = \left\{ \ba{cc} M^{-1}A^{T}M, & \mbox{bilinear form,}\\
               M^{-1}A^HM, & \mbox{ sesquilinear form.}
               \ea\right.\ee
Consider the Lie algebra $\L$ and the Jordan algebra $\J$ associated with the scalar product (\ref{sprod}) given by
\be \label{alg} \L :=\{ A\in \K^{n \times n}: A^{\star}=-A\} \mbox{ and } \J :=\{ A\in \K^{n \times n}: A^{\star}=A\}.\ee
In this paper, we consider $\S = \L $ or $\S= \J,$ and refer to the matrices in $\S$ as {\em structured matrices.}
The Jordan and Lie algebras so defined provide a general framework for analyzing a great many important classes of structured matrices including Hamiltonian, skew-Hamiltonian, symmetric, skew-symmetric, pseudosymmetric, persymmetric, Hermitian, Skew-Hermitian, pseudo-Hermitian, pseudo-skew-Hermitian matrices, to name only a few, see~(Table~2.1, \cite{fran6}).

For the rest of the paper, we set \beano \sym &:=&
\{A\in \K^{n \times n} : A^T = A\}, \, \,\, \sksym  \,\, := \,\, \{A\in \K^{n \times n} : A^T = -A\},  \\ \herm &:=& \{A\in \C^{n \times n} :
A^H = A\}, \, \,\, \skherm \,\, :=\,\,  \{A\in\C^{n \times n } : A^H = -A\}. \eeano
Also define the set $M\S$ by  $ M\S :=\{ M A : A \in \S\}.$ Then, in view of (\ref{adj}) and (\ref{alg}), it follows that
\be \label{str4} \S \in \{ \J, \L\} \Longleftrightarrow M \S \in \{ \sym, \, \sksym, \, \herm, \, \skherm\}.\ee
This shows that the four classes of structured matrices, namely, symmetric, skew-symmetric, Hermitian and skew-Hermitian matrices are prototypes of more general structured matrices belonging to the Jordan and Lie algebras given in (\ref{alg}).

Let $A, B, C$ and $D$ be matrices.  Then the matrix $T :=\bmatrix{ A & C\\
B & D}$ is called a {\bf dilation} of $A.$ The norm preserving
dilation problem is then stated as follows. Given matrices $A, B, C$
and a positive number $$ \mu \geq \max\left( \left\|\bmatrix{ A\\
B}\right\|_2, \, \left\|\bmatrix{ A & C}\right\|_2\right),
\eqno(1)$$ find all possible $D$ such that $\left\|\bmatrix{ A &
C\\ B & D}\right\|_2 \leq \mu.$ 

\begin{theorem}[Davis-Kahan-Weinberger, \cite{D}] \label{dkw}
Let $A, B, C$ be given matrices. Then for any positive number
$\mu$ satisfying (1),
  there exists $D$ such that $\left\|\bmatrix{  A & C \\
  B & D
}\right\|_2 \leq \mu.$ Indeed, those $D$ which have this property
are exactly those of the form $$D= - KA^HL +
\mu(I-KK^H)^{1/2}Z(I-L^HL)^{1/2},$$ where $K^H :=(\mu^2 \Ib -A^HA
)^{-1/2}B^H,~ L :=(\mu^2 \Ib -AA^H)^{-1/2}C$ and $Z$ is an
arbitrary contraction, that is, $\|Z\|_2 \leq 1.$
\end{theorem}

We mention that when $(\mu^2 I-A^HA)$ is singular, the
inverses in $K^H$ and $L$ are replaced by their Moore-Penrose
pseudo-inverses~(see, \cite{meinguet}). An interesting fact about Theorem~\ref{dkw} is that if $T(D)$ is symmetric or skew-symmetric or Hermitian or skew-Hermitian then the solution matrices $D$ are, respectively, symmetric or skew-symmetric or Hermitian or skew-Hermitian~\cite{Ba:thesis}.

\section{Solution of structured mapping problem }
For compact representations of our results, in the rest of the paper, we write $A^*$ to denote the transpose $A^T$ or the conjugate transpose $A^H.$  Often we write $A^*$ with $* \in \{T, H\}.$  With this notational convention, define the map $\F_* :
\K^{n \times p} \times \K^{n\times p}
\rightarrow \K^{n\times n}$ by \be \label{solmap} \F_*(X,B) := \left\{%
\begin{array}{ll}
    BX^\dagger + (BX^\dagger)^* - (X^\dagger)^* (X^*B)^*X^\dagger, \,\, \mbox{if}
    \,\,(X^*B)^*  =X^*B \\
    BX^\dagger - (BX^\dagger)^* -
(X^\dagger)^* (X^*B)^*X^\dagger, \,\, \mbox{if} \,\,(X^*B)^* =-X^*B  \\
BX^\dagger, \,\, \mbox{else} \\
\end{array}%
\right. \ee where $X^{\dagger}$ is the Moore-Penrose pseudoinverse
of $X$ and $* \in \{T, H\}.$

We write $\F_* = \F_{T}$ when $*=T$ and  $\F_* = \F_{H}$ when $*=H.$ Then it follows that $\F_*$ has the following properties. \begin{itemize}
\item[1.]  If  $X^TB$ is symmetric/skew-symmetric then $ \F_T(X, B) \in \sym/\sksym.$
\item[2.]  If  $X^HB$ is Hermitian/skew-Hermitian then $  \F_H(X, B) \in \herm/\skherm$.
\item[3.] If $BX^\dagger X = B$ then $ \F_*(X,B) X = B.$
\end{itemize}

This shows that  the matrix $\F_*(X, B)$ is a potential candidate
for solution of the structured mapping problem for four classes of structured matrices, namely, $ \sym, \sksym, \herm$ and $\skherm.$
More generally, the  following result provides a necessary and sufficient condition for
existence of solution of the structured mapping problem.


\begin{theorem}[Existence]\label{existlj}
Let $(X,B)\in\K^{n\times p}\times\K^{n \times p}$ and $\S \in \{
\J, \L\}.$  Then there is a matrix $A \in \S$ such that $AX = B$
if and only if (a) $BX^\dag X= B$ and (b) the condition in Table~\ref{extb} holds.
 \begin{table}[h]
\begin{center}\renewcommand{\arraystretch}{1.2}
\begin{tabular}{|c|c|c|}
  \hline
   $M$ & $\S =\J$ & $\S=\L$ \\
  \hline\hline
    $M^T=M$ &  $(X^T M B)^T=X^TMB$ &  $(X^TM B)^T=-X^TMB$ \\
  \hline
  $M^T=-M$ & $(X^TM B)^T=-X^TMB$ & $(X^TM B)^T=X^TMB$ \\
  \hline \hline
  $M^H=M$ & $(X^HM B)^H=X^HMB$ & $(X^HM B)^H=-X^HMB$ \\
  \hline
  $M^H=-M$ & $(X^HM B)^H=-X^HMB$ & $(X^HM B)^H=X^HMB$ \\
  \hline
\end{tabular} \caption{ \label{extb} Condition for $\S(X,B) \neq \emptyset.$ }
\end{center}
\end{table}
\end{theorem}

\noin\pf Suppose that there exists $A \in \S$ such that $AX =B.$ Then  $ X^*MB = X^*MAX$ for $ * \in \{T, H\}.$ Since $MA \in M\S,$ by (\ref{str4}) it follows that $ X^*MB$ is symmetric/skew-symmetric (resp., Hermitian/skew-Hermitian) when $ *=T$ (resp., $*=H$). Hence the conditions in (b) are satisfied. Again since $AX = B,$ we have  $BX^\dagger X = AXX^\dagger X = AX =B.$

Conversely, suppose that the conditions are satisfied. Then setting
$A := M^{-1}\F_*(X,MB),$ it follows from (\ref{str4}) and the properties of $\F_*$ that $ AX =B$ and $ A \in \S,$ where $*\in \{T, H \}.$  This completes the proof. $\blacksquare$

\vone
\begin{remark} \label{fullrank}
If $X$ has full column rank then $X^\dag X =I.$ Consequently, $BX^\dag X  = B.$  Thus for a full column rank matrix $X$ the condition (a) in Theorem~\ref{existlj} is automatically satisfied.
\end{remark}

\vone We mention that for the special case when $x \in \K^n$ and $b
\in \K^n,$ by Theorem~\ref{existlj},  we obtain the following
necessary and sufficient condition provided in~\cite{fran6}:

\begin{table}[h]
\begin{center}\renewcommand{\arraystretch}{1.2}
\begin{tabular}{|c|c|c|}
  \hline
  $M$ & $\S=\J$ & $\S=\L$ \\
  \hline\hline
  $M^T=M$ &  $\mbox{ any } x,b\in\K^n$ &  $x^TMb =0$ \\
  \hline
  $M^T=-M$ & $ x^TMb =0$  & $ \mbox{ any } x,b\in\K^n$  \\
  \hline \hline
  $M^H=M$ &  $x^HMb \in\R$ &  $x^HMb \in i \R$ \\
  \hline
 $M^H=-M$ &  $x^HMb \in i \R$ &  $x^HMb \in \R$ \\
   \hline
\end{tabular} \caption{\label{vec} Necessary and sufficient condition for $\S(x, b) \neq \emptyset$}
\end{center}
\end{table}

Now given a pair of matrices $X$ and $B$ in $\K^{n\times p}$ satisfying the conditions in Theorem~\ref{existlj},
the following result characterizes solution of the structured mapping problem.

\begin{theorem}[Characterization]\label{gsmp} Let
$(X,B)\in \K^{n\times p} \times \K^{n\times p}$ and $\S \in \{\J, \L\}.$ Suppose that
$\S (X,B) \neq \emptyset.$ Then  $A \in \S(X, B)$ if and only if $A$ is of
the form $$ A = \left\{%
\begin{array}{lll}
   M^{-1}\F_T(X, MB) + M^{-1}(I-XX^\dag)^T Z (I-
XX^\dag),  &\mbox{ if } M\S
\in\{\sym, \sksym \}, \\
  M^{-1}\G(X, MB) + M^{-1}(I-XX^\dag) Z (I-
XX^\dag), & \mbox{ if } M\S \in\{ \herm, \skherm\}, \\
\end{array}%
\right.$$ for some $Z$ such that $ Z \in M\S.$

 {\bf (a) Frobenius norm: } Define $ A_o :=M^{-1}\F_T(X, MB)$ when $M\S \in \{\sym, \sksym\},$ and  $A_o := M^{-1}\G(X, MB)$ when $M \S \in\{\herm, \skherm\}.$ Then $A_o$ is the unique matrix in $\S(X,B)$ such that $A_o X = B$ and that
$$ \sigma^\S(X, B) = \|A_o\|_F = \left\{%
\begin{array}{ll}
    \sqrt{2 \|BX^\dagger\|_F^2 - \tr(MBX^\dagger
(MBX^\dagger)^H (XX^\dagger)^T)}, \\  \hfill{\mbox{ when } M\S
\in\{\sym, \sksym \},} \\ & \\
  \sqrt{2 \|BX^\dagger\|_F^2 - \tr(MBX^\dagger
(MBX^\dagger)^H XX^\dagger)}, \\  \hfill{\mbox{ when } M\S
\in \{\herm, \skherm\}.} \\
\end{array}%
\right.$$

 {\bf (b) Spectral norm:}  Set $\sigma^{\rm unstruc}(X, B) := \sigma^\S(X, B)$ when $\S= \K^{n\times n}.$  Then   $$ \sigma^\S(X, B) = \|BX^\dagger\|_2 = \sigma^{\rm unstruc}(X, B).$$
 Suppose that $\rank(X) =r.$  Consider the SVD $X =U\Sigma V^H$ and partition $U$ as  $U=[U_1 \, \, U_2], $
where $U_1 \in \K^{n\times r}.$

\vone \noin{\bf Case-I.} If $M\S \in \{\sym, \sksym\}$ then consider the matrix
$$A_o :=M^{-1}\F_T(X,MB)-M^{-1}(I-XX^\dag)^TKU_1^H\overline{MBX^\dag
U_1}K^T(I-XX^\dag)+ M^{-1} f(Z),$$ where
$f(Z) =\mu\overline{U}_2(I-U_2^TKK^H\overline{U}_2)^{1/2}Z(I-U_2^H
\overline{K}K^TU_2)^{1/2}U_2^H$, $\mu=\|BX^\dagger\|_2,$ $$K=\left\{%
\begin{array}{ll}
   MBX^\dagger U_1 (\mu^2 I -
U_1^H\overline{MBX^\dagger}MBX^\dagger U_1)^{-1/2}, & \hbox{\mbox{when}\,\, $M\S = \sym,$} \\
   MBX^\dagger U_1 (\mu^2 I +
U_1^H\overline{MBX^\dagger}MBX^\dagger U_1)^{-1/2}, & \hbox{\mbox{when}\,\, $M\S = \sksym$,} \\
\end{array}%
\right.$$ and $Z$ is an arbitrary contraction such that $Z= Z^T$ (resp., $Z= -Z^T$) when $M\S = \sym$ (resp. $M\S = \sksym$). Then  $A_o \in \S(X, B)$ and $\|A_o\|_2 = \sigma^\S(X, B).$

\vone
\noin{\bf Case-II.}If $M \S = \herm$ then consider the matrix
$$A_o :=M^{-1}\G(X,MB)-M^{-1}(I-XX^\dag)KU_1^H MBX^\dag
U_1K^H(I-XX^\dag)+ M^{-1}f(Z),$$ where $ f(Z) = \mu U_2
(I-U_2^HKK^HU_2)^{1/2}Z(I-U_2^H KK^HU_2)^{1/2}U_2^H,$   $\mu=\|BX^\dagger\|_2,$
 $$K = MBX^\dagger U_1 (\mu^2 I -U_1^HMBX^\dagger MBX^\dagger U_1)^{-1/2},$$
and $Z = Z^H$ is an arbitrary contraction. Then  $A_o \in \S(X, B)$ and $\|A_o\|_2 = \sigma^\S(X, B).$\\

If $ M\S = \skherm$ then define $ A_o := -i \mathcal{G}(X, i B), $ where $\mathcal{G}(Y, D) \in \S(Y, D)$ is such that $\|\mathcal{G}(Y,D)\|_2 = \sigma^\S(Y, D)$ when $M\S = \herm.$ Then $ A_o \in \S(X, B)$ and $\|A_o\|_2 =\sigma^\S(X, B).$
\end{theorem}

\noin
\pf First, observe that $ AX = B \iff MAX = MB.$ Consequently, we have  $ A \in \S(X, B) \iff MA \in M\S(X, MB).$ Further, since $M$ is unitary,  $\sigma^\S(X, B) = \sigma^{M\S}(X, MB)$ for the spectral and the Frobenius norms. By (\ref{str4}), $ MS \in \{ \sym, \sksym, \herm, \skherm\}.$
Thus, it boils down to proving the results for symmetric, skew-symmetric, Hermitian and skew-Hermitian matrices.
Indeed, if $\phi(X, B)$ is a symmetric or skew-symmetric or Hermitian or skew-Hermitian solution of $AX =B$ then $ M^{-1}\phi(X, MB) \in \S(X, B).$

Further, note that a skew-Hermitian solution of $ AX = B$ can be obtained from a Hermitian solution of $ iAX = iB$ and vice versa. Indeed, if  $A X = B$ and $ X^HB$ is skew-Hermitian then $ i A X = iB$ and $ iX^HB$ is Hermitian. Hence $\phi(X, iB)$ is a Hermitian solution of $iAX =iB$ and  $ -i \phi(X, iB)$ is a skew-Hermitian solution of $ AX = B.$ Consequently, we only need to prove the results for symmetric, skew-symmetric and Hermitian matrices. We prove these results separately in Theorem~\ref{splsmp}. Hence the proof. $\blacksquare$

\vone \begin{remark} (a) We mention that the solution set $\S(X, B)$ as characterized in Theorem~\ref{gsmp} can be written compactly as
$${\small  \S(X, B) = \left\{%
\begin{array}{lll}
   M^{-1}\F_T(X, MB) + M^{-1}(I-XX^\dag)^T M\S (I-
XX^\dag),  &\mbox{ if } M\S
\in\{\sym, \sksym \}, \\
  M^{-1}\G(X, MB) + M^{-1}(I-XX^\dag) M\S (I-
XX^\dag), & \mbox{ if } M\S \in\{ \herm, \skherm\}. \\
\end{array}%
\right. }$$
Here $x+S := \{ x+ s : s \in S\}.$

\vone (b) We also mention that when $ X$ has full column rank, for the Frobenius norm, we have  .
$$ \sigma^{\S}(X, B) = \sqrt{2\|B(X^HX)^{-1/2}\|_F^2 -
\|(X^T\overline{X})^{-1/2}X^TMB(X^HX)^{-1/2}\|_F^2} $$ when $ M\S  \in \{\sym, \sksym\},$ and
$$\sigma^{\S}(X, B) = \sqrt{2\|B(X^HX)^{-1/2}\|_F^2 -
\|(X^HX)^{-1/2}X^HMB(X^HX)^{-1/2}\|_F^2}$$ when $ M\S \in \{ \herm, \skherm\}.$
Indeed, when $X$ has full column rank, we have $ X^\dagger X = I$ and  $ X^\dagger = (X^HX)^{-1} X^H$ and hence the results follow.
On the other hand, for the spectral norm, we have $\sigma^\S(X, B) = \|B(X^HX)^{-1/2}\|_2.$

\end{remark}

\vone
To complete the proof of Theorem~\ref{gsmp}, for the rest of this section, we consider $\S$ to be such that
 $\S \in \{\sym, \sksym, \herm, \skherm \}.$ Observe that if  $A \in
\K^{n \times n}$ is given by  \be\label{fronorm} A := \bmatrix{A_{11} & \pm A_{12}^* \\
A_{12} & A_{22}} \mbox{  then  } \|A\|_F = \left(2\left\|
\bmatrix{A_{11} \\ A_{12}} \right\|_F^2 - \|A_{11}\|_F^2 +
\|A_{22}\|_F^2\right)^{1/2}.\ee We
repeatedly use this fact in the sequel.

\begin{theorem}[Special solutions]\label{splsmp} Let
$(X,B)\in \K^{n\times p} \times \K^{n\times p}$ and $\S \in \{\sym, \sksym, \herm\}.$ Suppose that
$\S (X,B) \neq \emptyset.$ Then  $A \in \S(X, B)$ if and only if $A$ is of
the form $$ A = \left\{%
\begin{array}{lll}
   \F_T(X, B) + (I-XX^\dag)^T Z (I-
XX^\dag),  &\mbox{ if } \S
\in\{\sym, \sksym \}, \\
  \G(X, B) + (I-XX^\dag) Z (I-
XX^\dag), & \mbox{ if } \S = \herm, \\
\end{array}%
\right.$$ for some $ Z \in \S.$

\von
 {\bf (a) Frobenius norm: } Consider $ A_o :=\F_T(X, B)$ when $\S \in \{\sym, \sksym\},$ and  $A_o := \G(X, B)$ when $ \S =\herm.$ Then $A_o$ is a unique matrix in $\S(X, B)$ such that
$$ \sigma^\S(X, B) = \|A_o\|_F = \left\{%
\begin{array}{ll}
    \sqrt{2 \|BX^\dagger\|_F^2 - \tr(BX^\dagger
(BX^\dagger)^H (XX^\dagger)^T)}, \\  \hfill{\mbox{ when } \S
\in\{\sym, \sksym \},} \\ & \\
  \sqrt{2 \|BX^\dagger\|_F^2 - \tr(BX^\dagger
(BX^\dagger)^H XX^\dagger)}, \\  \hfill{\mbox{ when } \S
=\herm.} \\
\end{array}%
\right.$$

 {\bf (b) Spectral norm:}  We have  $$ \sigma^\S(X, B) = \|BX^\dagger\|_2 = \sigma^{\rm unstruc}(X, B).$$
 Suppose that $\rank(X) =r.$  Consider the SVD $X =U\Sigma V^H$ and partition $U$ as  $U=[U_1 \, \, U_2], $
where $U_1 \in \K^{n\times r}.$

\vone \noin{\bf Case-I.} If $\S \in \{\sym, \sksym\}$ then consider the matrix
$$A_o :=\F_T(X, B)- (I-XX^\dag)^TKU_1^H\overline{BX^\dag
U_1}K^T(I-XX^\dag)+  f(Z),$$ where
$f(Z) =\mu\overline{U}_2(I-U_2^TKK^H\overline{U}_2)^{1/2}Z(I-U_2^H
\overline{K}K^TU_2)^{1/2}U_2^H$, $\mu=\|BX^\dagger\|_2,$ $$K=\left\{%
\begin{array}{ll}
   BX^\dagger U_1 (\mu^2 I -
U_1^H\overline{BX^\dagger}BX^\dagger U_1)^{-1/2}, & \hbox{\mbox{when}\,\, $\S = \sym,$} \\
   BX^\dagger U_1 (\mu^2 I +
U_1^H\overline{BX^\dagger}BX^\dagger U_1)^{-1/2}, & \hbox{\mbox{when}\,\, $\S = \sksym$,} \\
\end{array}%
\right.$$ and $Z$ is an arbitrary contraction such that $Z= Z^T$ (resp., $Z= -Z^T$) when $\S = \sym$ (resp. $\S = \sksym$). Then  $A_o \in \S(X, B)$ and $\|A_o\|_2 = \sigma^\S(X, B).$

\vone
\noin{\bf Case-II.} If $ \S = \herm$ then consider the matrix
$$A_o :=\G(X, B)- (I-XX^\dag)KU_1^H BX^\dag
U_1K^H(I-XX^\dag)+ f(Z),$$ where $ f(Z) = \mu U_2
(I-U_2^HKK^HU_2)^{1/2}Z(I-U_2^H KK^HU_2)^{1/2}U_2^H,$   $\mu=\|BX^\dagger\|_2,$
 $$K = BX^\dagger U_1 (\mu^2 I -U_1^H BX^\dagger BX^\dagger U_1)^{-1/2},$$
and $Z = Z^H$ is an arbitrary contraction. Then  $A_o \in \S(X, B)$ and $\|A_o\|_2 = \sigma^\S(X, B).$
\end{theorem}

\noin\pf First, suppose that  $\S = \sym.$ By assumption there
exists $A \in \S$ be such that $A X = B.$ Note that $\range(X) =
\range(U_1).$ Thus representing $A$ relative to the decomposition
$$A : \range(X) \oplus \range(X)^\perp \rightarrow
\range(\overline{X}) \oplus \range (\overline{X})^\perp, $$ we
have $A = \overline{U} \,\,\overline{U}^H A UU^H.$ Set
$\widehat{A} = U^T A U = \bmatrix{
  A_{11} & A_{12}^T \\
  A_{12} & A_{22} \\
}.$  Then $\widehat{A} \in \S$  and  $\| A \| = \|
\widehat{A} \|$ for the spectral and the Frobenius norms. Set $\Sigma_{1} := \Sigma(1:r, 1:r).$ Let $
V = [ V_1, V_2]$ be a conformal partition of $V$ such that $ X =
U_1\Sigma_1V_1.$  Now $AX = B \Rightarrow
\overline{U}\widehat{A}U^H X =B.$ This gives $$ \bmatrix{A_{11} & A_{12}^T \\
A_{12} & A_{22}} \bmatrix{U_1^H \\ U_2^H}X = U^T B =
\bmatrix{U_1^T \\ U_2^T}B  \Rightarrow \bmatrix{A_{11} & A_{12}^T \\
A_{12}
& A_{22}} \bmatrix{\Sigma_1 V_1^H \\ 0} = \bmatrix{U_1^TB \\
U_2^TB}$$ and consequently  $ A_{11} \Sigma_1 V_1^H =U_1^T B$ and $
A_{12}\Sigma_1 V_1^H = U_2^T B.$

Therefore, we have $ A_{11} = U_1^TBV_1\Sigma_1^{-1}=U_1^TBX^\dag
U_1, \, A_{12} = U_2^T B V_1 \Sigma_1^{-1}. $ Notice that $A_{11}$
is symmetric if and only if $X^TB = B^TX$ and $BX^\dagger X = B.$
Indeed $X^TB=B^TX$ gives
$\overline{V}_1\Sigma_1U_1^TB=B^TU_1\Sigma_1V_1^H$ and thus
$B^TU_1=\overline{V}_1\Sigma_1U_1^TBV_1\Sigma_1^{-1}.$ Now \beano
(U_1^TBX^\dag U_1)^T &=& U_1^T(X^\dag)^TB^TU_1 =
U_1^T\overline{U}_1\Sigma_1^{-1} V_1^TB^TU_1 = \Sigma_1^{-1}
V_1^T\overline{V}_1\Sigma_1U_1^TBV_1\Sigma_1^{-1}\\ &=&
U_1^TBV_1\Sigma_1^{-1} = U_1^TBX^\dag U_1\eeano as desired.
Thus we have \begin{eqnarray} \label{eqn:sym} \widehat{A} = \bmatrix{U_1^TBX^\dag U_1 & (U_2^T B V_1 \Sigma_1^{-1})^T \\
U_2^T B V_1 \Sigma_1^{-1} & A_{22}}
\end{eqnarray}
 Then by~(\ref{fronorm}) we have $\|\widehat{A}\|_F^2 =  2
\|BX^\dagger\|_F^2 - \tr(BX^\dagger (BX^\dagger)^H
(XX^\dagger)^T)+ \|A_{22}\|_F^2.$  Hence, for the Frobenius norm, setting  $A_{22} =0$ we
obtain a unique matrix
$$ A = \overline{U}\bmatrix{U_1^TBX^\dag U_1 & (U_2^T B V_1 \Sigma_1^{-1})^T \\
U_2^T B V_1 \Sigma_1^{-1} & 0}  \bmatrix{U_1^H \\ U_2^H} =
BX^\dagger + (BX^\dagger)^T - (XX^\dagger)^T BX^\dagger = \F_{\tt{T}}(X,
B)$$ such that $A \in \S(X, B)$ and $ \|A\|_F = \sqrt{ 2
\|BX^\dagger\|_F^2 - \tr(BX^\dagger (BX^\dagger)^H
(XX^\dagger)^T)} = \sigma^\S(X, B).$

Now from (\ref{eqn:sym}) we have
\begin{eqnarray}\label{symm:strmap1} A &=&
\overline{U}\bmatrix{U_1^TBX^\dag U_1 & (U_2^T B V_1 \Sigma_1^{-1})^T \\
U_2^T B V_1 \Sigma_1^{-1} & A_{22}} \bmatrix{U_1^H \\ U_2^H} \nonumber\\
&=& \bmatrix{ \overline{U}_1 & \overline{U}_2}\bmatrix{U_1^TBX^\dag U_1U_1^H + \Sigma_1^{-1}V_1^T B^T U_2 U_2^H \\
U_2^T B V_1 \Sigma_1^{-1}U_1^H +  A_{22}U_2^H} \nonumber\\ &=&
(U_1U_1^H)^TBX^\dagger + (V_1\Sigma_1^{-1}U_1^H)^T B^T (I -
U_1U_1^H) + (I - U_1U_1^H)^TBX^\dagger + \overline{U}_2
A_{22}U_2^H \nonumber\\ &=& (XX^\dagger)^TBX^\dagger +
(X^\dagger)^TB^T(I - XX^\dagger) + (I - XX^\dagger)^TBX^\dagger +
\overline{U}_2A_{22}U_2^H  \nonumber\\ &=& BX^\dagger +
(BX^\dagger)^T - (XX^\dagger)^T BX^\dagger + \overline{U}_2 U_2^T
\overline{U}_2 A_{22} U_2^H U_2 U_2^H \nonumber\\ &=& BX^\dagger +
(BX^\dagger)^T - (X^\dagger)^T X^TBX^\dagger + (I - XX^\dagger)^T
Z (I - XX^\dagger) \nonumber \\ &=& \F_{\tt{T}}(X, B)+ (I - XX^\dagger)^T Z
(I - XX^\dagger),\nonumber \end{eqnarray} where $Z \in \S$ is
arbitrary.

For the spectral norm, again consider the matrix $\widehat{A}$ given in (\ref{eqn:sym}) and set $$ \mu
 := \left\|\bmatrix{U_1^TBV_1\Sigma_1^{-1} \\
U_2^T B V_1 \Sigma_1^{-1} }\right\|_2 = \| U^T
BV_1\Sigma_1^{-1}\|_2 =\|BX^\dagger \|_2. $$ Then it follows that
$\|\widehat{A}\|_2 \geq \mu.$ Now by
Theorem~\ref{dkw} we have $\|\widehat{A}\|_2 = \mu$ when
\be \label{dl} A_{22}= -K_1\overline{A}_{11}K_1^T + \mu
(I-K_1K_1^H)^{1/2}Z(I-\overline{K}_1K_1^T)^{1/2} \ee where $ K_1
=U_2^TBX^\dag U_1 (\mu^2 I - U_1^H\overline{BX^\dagger} BX^\dag
U_1)^{-1/2}$ and $Z$ is an arbitrary contraction such that $ Z = Z^T$ (resp., $Z= -Z^T$)  when $\S= \sym$ (resp., $\S=\sksym$). Thus we
have $$A= BX^\dagger + (BX^\dagger)^T - (XX^\dagger)^T BX^\dagger +
\overline{U}_2A_{22}U_2^H $$ such that $ A  \in \S,$ $ AX = B$ and
$\|A\|_2 = \|BX^\dag\|_2 = \sigma^\S(X, B),$ where $A_{22}$ is
given in (\ref{dl}).  Upon simplification we obtain the
desired form of $A.$

\vone

Next suppose that  $\S = \sksym.$ Again representing $A$ relative to the decomposition
$$A : \range(X) \oplus \range(X)^\perp \rightarrow
\range(\overline{X}) \oplus \range (\overline{X})^\perp, $$ we
have $A = \overline{U} \,\,\overline{U}^H A UU^H.$ Set
$\widehat{A} = U^T A U = \bmatrix{
  A_{11} & -A_{12}^T \\
  A_{12} & A_{22} \\
}.$  Then $\widehat{A} \in \S$  and  $\| A \| = \|
\widehat{A} \|$ for the spectral and the Frobenius norms.
By assumption  $AX = B \Rightarrow
\overline{U}\widehat{A}U^H X =B.$ This gives $$ \bmatrix{A_{11} & -A_{12}^T \\
A_{12} & A_{22}} \bmatrix{U_1^H \\ U_2^H}X = U^T B =
\bmatrix{U_1^T \\ U_2^T}B  \Rightarrow \bmatrix{A_{11} & - A_{12}^T \\
A_{12}
& A_{22}} \bmatrix{\Sigma_1 V_1^H \\ 0} = \bmatrix{U_1^TB \\
U_2^TB}.$$ Consequently we have  $ A_{11} \Sigma_1 V_1^H =U_1^T B$ and $
A_{12}\Sigma_1 V_1^H = U_2^T B$ and hence   \begin{eqnarray} \label{eqn:sksym} \widehat{A} =
 \bmatrix{U_1^TBX^\dag U_1 & -(U_2^T B V_1 \Sigma_1^{-1})^T \\ U_2^T B V_1 \Sigma_1^{-1} & A_{22}}.
\end{eqnarray}  As before, setting $A_{22} =0$ in (\ref{eqn:sksym}) we obtain the desired results for the Frobenius norm.

As for the spectral norm, applying Theorem~\ref{dkw} to the matrix $\widehat{A}$ in (\ref{eqn:sksym}) and following steps similar to those in the case when $\S = \sym,$ we obtain the desired results.

\vone
 Finally, suppose that  $\S = \herm.$ Then  representing $A$ relative to the decomposition
$$A : \range(X) \oplus \range(X)^\perp \rightarrow
\range({X}) \oplus \range ({X})^\perp, $$ we have $A = {U}
\,\,{U}^H A UU^H.$ Set $\widehat{A} = U^H A U = \bmatrix{
  A_{11} & A_{12}^H \\
  A_{12} & A_{22} \\
}.$  Then $\widehat{A} \in \S$  and  $\| A \| = \|
\widehat{A} \|$ for the spectral and the Frobenius norms.
Now $AX = B \Rightarrow
{U}\widehat{A}U^H X =B.$ This gives $$ \bmatrix{A_{11} & A_{12}^H \\
A_{12} & A_{22}} \bmatrix{U_1^H \\ U_2^H}X = U^H B =
\bmatrix{U_1^H \\ U_2^H}B  \Rightarrow \bmatrix{A_{11} & A_{12}^H \\
A_{12}
& A_{22}} \bmatrix{\Sigma_1 V_1^H \\ 0} = \bmatrix{U_1^HB \\
U_2^HB}.$$ Consequently, we have  $ A_{11} \Sigma_1 V_1^H =U_1^H B$ and $
A_{12}\Sigma_1 V_1^H = U_2^H B$ and hence
 \begin{eqnarray} \label{eqn:herm} \widehat{A} = \bmatrix{U_1^HBX^\dag U_1 & (U_2^H B V_1 \Sigma_1^{-1})^H \\
U_2^H B V_1 \Sigma_1^{-1} & A_{22}}
\end{eqnarray}
 Now by~(\ref{fronorm}) we have $\|\widehat{A}\|_F^2 =  2
\|BX^\dagger\|_F^2 - \tr(BX^\dagger (BX^\dagger)^H XX^\dagger)+
\|A_{22}\|_F^2.$  Hence setting  $A_{22} =0$ in (\ref{eqn:herm}), we obtain a unique
matrix
$$ A = {U}\bmatrix{U_1^HBX^\dag U_1 & (U_2^H B V_1 \Sigma_1^{-1})^H \\
U_2^H B V_1 \Sigma_1^{-1} & 0}  \bmatrix{U_1^H \\ U_2^H} =
BX^\dagger + (BX^\dagger)^H - XX^\dagger BX^\dagger = \G(X, B)$$
such that $A \in \S(X, B)$  and $ \|A\|_F = \sqrt{ 2
\|BX^\dagger\|_F^2 - \tr(BX^\dagger (BX^\dagger)^H XX^\dagger)} =
\sigma^\S(X, B).$

Now from (\ref{eqn:herm}) we have
\begin{eqnarray}\label{symm:strmap1} A &=&
{U}\bmatrix{U_1^HBX^\dag U_1 & (U_2^H B V_1 \Sigma_1^{-1})^H \\
U_2^H B V_1 \Sigma_1^{-1} & A_{22}} \bmatrix{U_1^H \\ U_2^H} \nonumber\\
&=& \bmatrix{ {U}_1 & {U}_2}\bmatrix{U_1^HBX^\dag U_1U_1^H + \Sigma_1^{-1}V_1^H B^H U_2 U_2^H \\
U_2^H B V_1 \Sigma_1^{-1}U_1^H +  A_{22}U_2^H} \nonumber\\ &=&
(U_1U_1^H)BX^\dagger + (V_1\Sigma_1^{-1}U_1^H)^H B^H (I -
U_1U_1^H) + (I - U_1U_1^H)BX^\dagger + {U}_2 A_{22}U_2^H
\nonumber\\ &=& XX^\dagger BX^\dagger + (X^\dagger)^HB^H(I -
XX^\dagger) + (I - XX^\dagger)BX^\dagger + {U}_2A_{22}U_2^H
\nonumber\\  &=& BX^\dagger + (BX^\dagger)^H - (X^\dagger)^H
X^HBX^\dagger + (I - XX^\dagger) Z (I - XX^\dagger) \nonumber \\
&=& \G(X, B)+ (I - XX^\dagger) Z (I - XX^\dagger),\nonumber
\end{eqnarray} where $Z \in \S$ is arbitrary.

For the spectral norm, again consider the matrix $\widehat{A}$ given in(\ref{eqn:herm}) and set
$$ \mu := \left\|\bmatrix{U_1^HBV_1\Sigma_1^{-1} \\
U_2^H B V_1 \Sigma_1^{-1} }\right\|_2 = \| U^H
BV_1\Sigma_1^{-1}\|_2 =\|BX^\dagger \|_2. $$ Then it follows that
$\|\widehat{A}\|_2 \geq \mu.$ Now by
Theorem~\ref{dkw} we have $\|\widehat{A}\|_2 = \mu$ when
\be \label{dlherm} A_{22}= -K_1{A}_{11}K_1^H + \mu
(I-K_1K_1^H)^{1/2}Z(I-{K}_1K_1^H)^{1/2}\ee where $ K_1
=U_2^HBX^\dag U_1 (\mu^2 I - U_1^H{BX^\dagger} BX^\dag
U_1)^{-1/2}$ and $Z =Z^H$ is an arbitrary contraction. Thus we
have  $A= BX^\dagger + (BX^\dagger)^H - XX^\dagger BX^\dagger +
U_2A_{22}U_2^H = \G(X, B) +U_2A_{22}U_2^H $ such that $ A \in \S(X, B)$
 and $\|A\|_2 = \|BX^\dag\|_2 = \sigma^\S(X, B),$ where
$A_{22}$ is given in (\ref{dlherm}). Finally, upon simplification, we obtain  the desired form of $A.$ This completes the proof. $\blacksquare$

\section{Applications of structured mapping problem}
We now consider a few applications of the structured mapping problem. As before, we consider $ \S \in \{\J, \L\}.$ Given a full column rank matrix $ X \in \K^{n\times p},$  a matrix $ D \in \K^{p\times p}$ and a structured matrix $A \in \S,$ we say that $ (X, D)$ is an {\em invariant pair} for $A$ if $ AX = XD.$
Then a partially prescribed inverse eigenvalue problem can be stated as follows.

\vone \noin {\sc Problem-I.}  Given a full column rank matrix $ X \in \K^{n\times p}$ and a matrix $ D \in \K^{p\times p},$ find a matrix $A_o \in \S,$ if it exists,
such that $ A_oX = XD$ and $\|A_o\| = \tau^\S(X, D),$ where  $$\tau^\S(X, D) := \inf\{ \|A\|: A \in \S \mbox{ and } AX = X D\}.$$

When  $A_o \in \S$ exists, $(X, D)$ provides a partial spectral decomposition of $A_o$ in the sense that  $\Lambda(D) \subset \Lambda(A)$ and that $\mathrm{span}(X)$ is an invariant subspace of $A$ corresponding to $\Lambda(D).$ Obviously, {\sc Problem-I} is a structured mapping problem for $X$ and $ B := XD.$ Consequently, it has a solution if and only if $X$ and $B:= XD$ satisfy the conditions in Theorem~\ref{existlj}, that is, if and only if $(X^*MXD)^* = X^*MXD$ or $(X^*MXD)^* = -X^*MXD$ for $* \in \{T, H\}.$  An optimal solution $ A_o$ is then given by Theorem~\ref{gsmp} with $ \tau^\S(X, D)= \|A_o\|= \sigma^{\S}(X, XD)$ for the spectral and the Frobenius norms. Note that {\sc Problem-I} has a solution when $ X^*MX =0$ and in such a case the subspace spanned by the columns of $X$ is called $M$-neutral - in short $X$ is $M$-neutral. Thus if $X$ is $M$-neutral then for any $D \in \K^{p\times p},$ $(X, D)$ is an invariant pair for some $A \in \S.$

A related problem which arises when analyzing backward errors of approximate invariant pairs is as follows.

\vone \noin {\sc Problem-II.}  Given a full column rank matrix $ X \in \K^{n\times p}, $ a matrix $ D \in \K^{p\times p}$ and a structured matrix $A \in \S,$ find a structured matrix $E_o \in \S,$ if it exists,  such that $ E_oX = XD$ and $ \|E_o-A\| = \eta^\S(A, X,D),$ where
  $ \eta^\S(A, X,D) :=\inf\{ \|E-A\|: E \in \S \mbox{ and } EX = X D\}.$

\vone
{\sc Problem-II} occurs, for example, when Krylov subspace method such as the Arnoldi method is used to compute a few eigenvalues of a (large)  matrix $A.$  Indeed, starting with a unit vector $v_1 \in \C^n,$ after $p$ steps of Arnoldi method we have $$ AV_p = V_pH_p+ \alpha v_{p+1} e_{p+1}^T,$$ where $V_p :=[v_1, \ldots, v_p],$ $H_p$ is upper hessenberg and $\{v_1, \ldots, v_{p+1}\}$ is an orthonormal basis of the Krylov subspace $ \mathcal{K}(v_1, A) := \mathrm{span}( v_1, Av_1, \ldots, A^pv_1).$ Thus, when $|\alpha|$ is small, $(V_p, H_p)$ is an approximate invariant pair of $A$ and hence the backward error of $(V_p, H_p)$ may be gainfully used in analyzing errors in the computed eigenvalues of $A.$

Writing $ E = A +\Delta A,$ it follows that {\sc Problem-II} is a structured mapping problem for $X$ and $ B:= XD -AX.$ Indeed, if $ \Delta A \in \S$ is such that $ (A+\Delta A) X = X D$ then $ \Delta A X = XD-AX =B.$ Consequently, $\Delta A \in \S$ exists if and only if $X$ and $B= XD -AX$ satisfy the conditions in Theorem~\ref{existlj}.  An optimal solution $ \Delta A_o$ is then given by Theorem~\ref{gsmp} with $$\eta^\S(A, X, D) = \|\Delta A_o\| = \sigma^{\S}(X, XD-AX)$$ for the spectral and the Frobenius norms. The quantity $\eta^\S(A, X, D)$ is the {\em structured backward error} of $(X, D)$ as an approximate invariant pair of $A.$
Note that the conditions in Theorem~\ref{existlj} are satisfied if and only if $ X^HMXD -X^HMAX$ (resp., $ X^TMXD -X^TMAX$) is Hermitian or skew-Hermitian (resp., symmetric or skew-symmetric). In particular, this condition is satisfied when $ M^H= -M,$  $ \S = \L$  and $ X \in \K^{n\times p}$ is $M$-neutral. This fact plays an important role in spectral perturbation analysis of Hamiltonian matrices, see~\cite{AlaBKMM}.

We now consider a special case when $X = x \in \C^n$ and $D = \lam \in \C$ and derive structured backward error $\eta^\S(A, x, \lam)$ of an approximate eigenpair $(\lam, x)$ of $A.$  Set $ r : = \lam x - Ax.$ Then there is a matrix $E \in \S$ such
 that $(A +E)x = \lam x$ if and only if $ x$ and $r$ satisfy the
 condition in Table~\ref{vec}.  Let  $\eta^\S_2(A, x, \lam)$ (resp., $\eta^\S_F(A, x, \lam)$) denote $\eta^\S(A, x,\lam)$ for the spectral norm (resp., Frobenius norm).
  Then by Theorem~\ref{gsmp}, we have the following.
\begin{corollary}
Let $\S\in\{\J,\L\}$ and $A\in\S.$ Suppose that $\lam \in \C$ and $x \in \C^n\setminus\{0\}$ satisfy the
 condition in Table~\ref{vec}. Then
we have $$\eta^\S_2(A, x, \lam) =
 \|r\|_2/{\|x\|_2}\,\, \mbox{ and } \,\,\eta^\S_F(A, x, \lam) = \sqrt{ 2 \|r\|_2^2 -
 |\inp{r}{x}_M|^2}.$$ Define $E$ by
 $$ E := \left\{%
\begin{array}{ll}
    (x^Tr) M^{-1}\overline{x}x^H +  M^{-1} \overline{x}r^T(I-xx^H)
+ M^{-1}(I-\overline{x}x^T)rx^H, & \hbox{ if $MA\in \sym$ } \\
    M^{-1}(I-\overline{x}x^T)rx^H- M^{-1}
\overline{x}r^T(I-xx^H) , & \hbox{if $MA\in \sksym$} \\
    (x^Hr) M^{-1} xx^H + M^{-1} xr^H(I-xx^H) + M^{-1}(I-xx^H)rx^H, & \hbox{if $MA\in \herm$} \\
  \end{array}%
\right.$$ Then  $ E \in \S$ is a unique matrix  such that $(A+E)x = \lam x$ and
$\|E\|_F = \eta^\S_F(A, x, \lam).$ Further, when $\|r\|_2\neq |x^Hr|,$ define
 $$ \Delta A := \left\{%
\begin{array}{ll}
    E-\dfrac{\overline{x^Tr}
M^{-1}(I-xx^H)^Trr^T(I-xx^H)}{\|r\|_2^2-|x^Hr|^2}, & \hbox{if $MA\in \sym$} \\
    E, & \hbox{if $MA\in \sksym$} \\
    E-\dfrac{\overline{x^Hr}
M^{-1}(I-xx^H)rr^H(I-xx^H)}{\|r\|_2^2-|x^Hr|^2}, & \hbox{if $MA\in \herm$} \\
    \end{array}%
\right. $$ else set $ \Delta A := E.$ Then  $ \Delta A \in \S$ such that $(A+\Delta A)x =
\lam x$ and $\|\Delta A\|_2 = \eta^\S_2(A, x, \lam).$

\end{corollary}

\vone

Finally, yet another related problem which arises when dealing with backward errors of approximate invariant subspaces as well as in inverse eigenvalue problem with a specified invariant subspace is as follows.

\vone \noin  {\sc  Problem-III.}  Let $\mathcal{X}$ be a $p$-dimensional subspace of $\K^n$ and $A \in \S.$ Then find a structured matrix $E \in \S,$ if it exists,  such that $ E \mathcal{X} \subset \mathcal{X}$ and $ \|A-E\| = \omega^\S(A, \cal{X}),$ where
 $$ \omega^\S(A, \mathcal{X}) :=\inf\{ \|E-A\|: E \in \S \mbox{ and } E\mathcal{X} \subset\mathcal{X}\}.$$

Let $ U \in \K^{n\times p}$  be an isometry such that $\mathrm{span}(U) = \mathcal{X}.$ If $E \in \S$ is such that $ E\mathcal{X} \subset \mathcal{X}$ then  $ E U = U D$ for some $p$-by-$p$ matrix $D.$ Then setting  $\Delta A := E-A,$ we have $ \Delta A U = UD-AU$ and hence by Theorem~\ref{gsmp},  $ \omega^\S(A, \mathcal{X}) = \|AU-UD\|_2$ for the spectral norm. The choice of $D$ that minimizes $\|AU-UD\|$ is given by the following result.

\begin{proposition}\label{minnorm} Let  $ U \in \K^{n\times p}$ be an isometry and $ A \in \K^{n\times n}$. Set $ P := UU^H.$ Then for the spectral and the Frobenius norms, we have  $$ \min_{ D}\|AU-UD\| =\|AU - U (U^HAU)\| = \|(I-P) A P\|,$$ where the minimum is taken over $ \K^{p\times p}.$ Further, if $A$ is Hermitian (resp., skew-Hermitian) then so is the minimizer $D.$ On the other hand, if $A$ is symmetric (resp., skew-symmetric) then so is the minimizer $D$ provided $U$ is real.
\end{proposition}

\noin
\pf Let $ Z := [U, V] $ be unitary. Then for the spectral and the Frobenius norms, we have   $$\|AU -UD\| = \|Z^H(AU-UD)\| = \| \bmatrix{ U^HAU - D\\ V^HAU }\| \geq \|V^HAU\|$$ and the equality holds for $D= U^HAU.$  Now $\|V^HAU\| = \|VV^HAUU^HU\| = \|(I-P)AP\|$ gives the desired minimum. That $D$ inherits the structure of $A$ when $A$ is Hermitian or skew-Hermitian, and under the additional assumption of $U$ being real when $A$ is symmetric or skew-symmetric, is immediate.  $\blacksquare.$

\vone
Thus, for a solution of {\sc Problem~III}, we choose an isometry $U$ whose columns form an orthonormal basis of $\mathcal{X}$ and set $ D := U^HAU.$ Then {\sc Problem-III} reduces to {\sc Problem-II} for the pair $(U, D).$ However, since $U$ is an isometry and $D = U^HAU,$ existence criterion as well as solutions of {\sc Problem~III} have simpler forms. Indeed, when $M\S = \herm $ (resp., $M\S = \skherm$), we have $ U^HM(I-P)AU $ is Hermitian (resp., skew-Hermitian) and hence by Theorem~\ref{existlj},  {\sc Problem-III} has a solution. On the other hand, when $M\S = \sym $ (resp., $M\S = \sksym$), we have $U^T M(I-P)AU)$ is symmetric (resp., skew-symmetric) provided $U$ is real. Hence by Theorem~\ref{existlj},  {\sc Problem-III} has a solution. We mention once again that {\sc Problem-III} has a solution for the case when $M\S \in \{\sym, \sksym\}$ provided that the orthogonal projection $P$  is real. Thus, in either case, if $ E\in \S$ and $E\mathcal{X} \subset \mathcal{X}$ then by Theorem~\ref{gsmp}, we have $$ E = \left\{
                                        \begin{array}{ll}
                                          A+M^{-1}\mathcal{F}_H(U, M(I-P)AU) + M^{-1} (I-P) Z (I-P), & \mbox{ if } M\S \in \{\herm, \skherm\} \\
                                          A+M^{-1}\mathcal{F}_T(U,  M(I-P)AU) + M^{-1} (I-P) Z (I-P), & \mbox{ if } M\S \in \{\sym, \sksym\} \\
                                        \end{array}
                                      \right.$$
where $Z \in M\S$ is arbitrary. Further, we have the following result from Theorem~\ref{gsmp}.

\begin{theorem}\label{prob3} Let $ \mathcal{X}$ be a $p$-dimensional subspace of $\K^n$ and $ A \in \S.$ Let $P$ be the orthogonal projection on $\mathcal{X}$ given by $P:=UU^*.$
Define $ E_o := A+ M^{-1}\mathcal{F}_H(U, M(I-P)AU)$ when $  M\S \in \{\herm, \skherm\},$ and $ E_o := A+ M^{-1}\mathcal{F}_T(U,  M(I-P)AU)$ when $M\S \in \{\sym, \sksym\} $ and $U$ is real. Then $E_o \mathcal{X} \subset \mathcal{X}.$  Further, $ \omega^\S(A, \mathcal{X}) = \|A-E_o\|_2 = \|(I-P)AP\|_2$ for the spectral norm, and $$\omega^\S(A, \mathcal{X}) =\|A-E_o\|_F= \sqrt{2\|(I-P)AP\|_F^2-\|PM(I-P)AP\|_F^2}$$ for the Frobenius norm.

In particular, when $\S \in \{\herm, \skherm, \sym, \sksym\}$ we have
$E_o = PAP+ (I-P)A(I-P)$ and $ \omega^\S(A, \mathcal{X}) = \|A-E_o\|_2 = \|(I-P)AP\|_2$ for the spectral norm, and  $ \omega^\S(A, \mathcal{X}) = \|A-E_o\|_F = \sqrt{2}\, \|(I-P)AP\|_F$ for the Frobenius norm.
\end{theorem}

\vone
We mention that $E_o$ is a unique solution of {\sc Problem-III} for the Frobenius. By contrast, {\sc Problem-III} has infinitely many solutions for the spectral norm, which are of the form $E_o + G$ for appropriate $G \in \S$ as given in Theorem~\ref{gsmp}. We also mention that for the special case when, for example, $\S = \herm,$ the results in Theorem~\ref{prob3} can be obtained easily without resorting to Theorem~\ref{gsmp}. Indeed, if $E \in \herm$ and $E\mathcal{X} \subset \mathcal{X}$ then $ (I-P)EP=0 = PE(I-P)$ which gives $E = PEP+(I-P)E(I-P)$ or in (operator) matrix notation, we have  $$ E = \bmatrix{ PEP & 0 \\ 0 & (I-P)E(I-P)}.$$ Since $ A = PAP+ PA(I-P) + (I-P)AP+ (I-P)A(I-P)$ or in matrix notation $$A = \bmatrix{ PAP & PA(I-P) \\ (I-P)AP & (I-P)A(I-P)},$$ it follows that $ \|E-A\|$ is minimized for the spectral norm as well as for the Frobenius norm when $ E = PAP+ (I-P)A(I-P)$ and the minimum is given by $\|(I-P)AP\|_2$ for the spectral norm, and $\sqrt{2}\, \|(I-P)AP\|_F$ for the Frobenius norm.

%
%

\subsection{Structured pseudospectra} Let $A \in \S.$ Recall that $\eta^\S(A, x, \lam)$ is the structured backward error of $(\lam, x) \in \C\times \C^n$ as an approximate eigenpair of $A$ with the convention $\eta^\S(A,x, \lam) = \infty $ when $x$ and $ r:= \lam x- Ax$ do not satisfy the condition for structured mapping. For the spectral norm, set $$ \eta^\S(\lam, A) := \inf\{ \eta^\S(A, x, \lam) : x \in \C^n \mbox{ and }\|x\|_2=1\}.$$ Then $\eta^\S(\lam, A)$ is the structured backward error of $\lam$ as an approximate eigenvalue of $A.$ The backward error $\eta^\S(\lam, A)$ can be used to define structured pseudospectrum $\slam(A)$ of $A:$ $$ \slam(A) :=
\bigcup_{\|\Delta{A}\|_2 \leq \ep} \{ \Lambda(A+\Delta{A}) : \Delta{A}
\in \S\} = \{ \lam \in \C : \eta^\S(\lam, A) \leq \ep\}.$$ See~\cite{lntbook} for more on pseudospectra and \cite{graillat:pseudo, karow, rump} for structured pseudospectra.

We denote $\eta^\S(\lam, A)$ by $\eta(\lam, A)$ when $\S = \C^{n\times n}.$ Then obviously $\eta(\lam, A)$ is the unstructured backward error of $\lam$ as an approximate eigenvalue of $A$ and we have  $\eta(\lam, A) = \sigmin(A- \lam I),$ where $\sigmin(A)$ is the smallest singular value of $A.$   In contrast, for many important structures determination of $\eta^{\S}(\lam, A) $  is a hard optimization problem and hence computation of $\slam(A)$ is a challenging task~\cite{karow}. Nevertheless,  for certain structures it turns out that  $\eta^\S(\lam, A) = \eta(\lam, A)$ holds for all $\lam \in \C,$ and for some other structures the equality holds only for certain $\lam \in \C.$ We denote $\slam(A)$ by $\elam(A),$ the unstructured pseudospectrum of $A,$ when $ \S = \C^{n \times n}.$

\begin{theorem} \label{str1}
Let $\S := \J$ when $M^T =
M,$ and $ \S \in \{ \J,\L
\}$ when $M^T=-M.$ Let $ A \in \S.$ Then for $\lam \in \C,$ there
exists $\Delta{A} \in \S$ such that $ \lam \in \Lambda(A+\Delta{A})$ and
$\eta^{\S}(\lam, A) = \|\Delta{A}\|_2 = \eta(\lam, A).$ Consequently, we
have $$ \elam^{\S}(A) = \elam(A).$$
\end{theorem}

\vone\noin {\bf Proof:} Consider the case $\S := \J$ when $M^T =
M.$ Then  $MA \in \sym.$ Let $ \lam \in \C.$ It follows that $A-\lam I \in \S,$
that is, $(M(A-\lam I))^T = M(A-\lam I).$  Since $M(A-\lam I)$ is
complex symmetric there exists a unitary matrix $U$ such that the
symmetric Takagi factorization~\cite{horn:book} $ M(A-\lam I) =
U\Sigma U^T$ holds, where $\Sigma$ is a diagonal matrix containing
singular values of $M(A-\lam I)$ ordered in descending order of
magnitude. Note that $\eta(\lam, A)= \sigmin(A-\lam I) =
\sigmin(M(A-\lam I)) = \Sigma(n,n).$ Let $ u:= U(:, n).$ Then we
have $ M(A-\lam I) \overline{u} = \eta(\lam, A)\overline{u}.$ This
gives $ (A -\eta(\lam, A)M^{-1}uu^T)\overline{u} = \lam
\overline{u}.$
 Setting $\Delta{A} := -\eta(\lam, A)M^{-1}uu^T,$ we have $
 \Delta{A} \in \S$ and $\|\Delta{A}\|_2 = \eta(\lam, A) =
 \eta^{\S}(\lam, A)$. Hence the results follow.

Next, consider the case  $ \S \in \{ \J,\L
\}$ when $M^T=-M.$ Then $MA \in \{\sksym, \sym\}.$ First consider $MA\in\sksym.$ For
$\lam \in \C, $ we have $M(A-\lam I) \in \S,$ that is, $(M(A-\lam
I))^T = -M(A-\lam I).$ Since $M(A-\lam I) $ is complex
skew-symmetric, we have the skew-symmetric Takagi
factorization~\cite{horn:book}
$$ M(A-\lam I) = U\diag(d_1, \cdots, d_m) U^T,$$ where $U$ is
unitary, $d_j := \bmatrix{ 0 & s_j\\ -s_j & 0},$ $s_j \in \C$ is
nonzero and $|s_j|$ are singular values of $M(A-\lam I).$ Here the
blocks $d_j$ appear in descending order of magnitude of $|s_j|.$
Note that $M(A-\lam I) \overline{U} = U \diag(d_1, \cdots, d_m).$
Let $ u := U(:, n-1:n).$ Then $M(A-\lam I) \overline{u} = u d_m =
u d_m u^T\overline{u}.$ This gives $ (MA -ud_mu^T)\overline{u} =
\lam M \overline{u}.$ Hence taking $\Delta{A} := -M^{-1} ud_m
u^T,$ we have $ \lam \in \Lambda(A+\Delta{A}),\,\, \Delta{A} \in \S$
and $\|\Delta{A}\|_2 = |s_m| = \sigmin(M(A-\lam I)) = \sigmin(A-\lam
I) = \eta(\lam, A).$ Hence $\eta^{\S}(\lam, A) = \eta(\lam, A)$
and the desired result follows.

Finally, consider the case $M^T = -M$ and $ MA \in \sym.$ Let
$\lam \in \C.$ Note that $ \sigmin(M(A-\lam I)) = \sigmin(A-\lam
I) = \eta(\lam, A).$ Let $u$ and $v$ be unit left and right
singular vector of $M(A-\lam I)$ corresponding to $\eta(\lam, A).$
Then $M(A-\lam I)v =\eta(\lam, A) u.$ This gives $ Av-\eta(\lam,
A) M^{-1}u = \lam v.$ Let $ E \in \C^{n\times n}$ be such that $E=
E^T,$ $Ev = u$ and $\|E\| = 1.$ Such a matrix always exists (Theorem~\ref{gsmp}). Then
setting $\Delta{A} := -\eta(\lam, A) M^{-1} E,$ we have $
(A+\Delta{A})v = \lam v,$ $\Delta{A} \in \S$ and $\|\Delta{A}\|_2 =
\eta(\lam, A) = \eta^{\S}(\lam, A).$ Hence the result follows.
$\blacksquare$

\vone For Lie and Jordan algebras corresponding to sesquilinear
form induced by $M,$ we have partial equality between structured
and unstructured pseudospectra.

\begin{theorem} \label{str2} Let $\S := \J$ when $M^H =\pm M.$ Let $ A
\in \S.$ Then for $\lam \in \R,$ there exists $\Delta{A} \in \S$ such
that $ \lam \in \Lambda(A+\Delta{A})$ and $\eta^{\S}(\lam, A) =
\|\Delta{A}\|_2 = \eta(\lam, A).$ Consequently, we have $$
\elam^{\S}(A)\cap \R = \elam(A) \cap \R.$$

\noin Next, consider $\S :=\L$ when $M^H =\pm M.$  Let $ A \in \S.$ Then for $\lam
\in i\R,$ there exists $\Delta{A} \in \S$ such that $ \lam \in
\Lambda(A+\Delta{A})$ and $\eta^{\S}(\lam, A) = \|\Delta{A}\| = \eta(\lam,
A).$ Consequently, we have $$ \elam^{\S}(A)\cap i\R = \elam(A)
\cap i\R.$$

\end{theorem}

\noin{\bf Proof:} Consider the case $\S := \J$ when $M^H =\pm M.$ Then  $ MA\in\herm$ when $M= M^H$ and  $
MA\in\skherm$ and $M^H = -M.$ First consider $ MA\in\herm$ when $M= M^H.$ Now for $\lam \in \R,$  $M(A-\lam I) \in \S.$ Since $M(A-\lam I)$
is Hermitian, we have the spectral decomposition $M(A-\lam I) =
U\diag(\mu_1, \cdots, \mu_n)U^H,$ where $U$ is unitary and
$\mu_j$'s appear in descending order of their magnitudes. Note
that $|\mu_n| = \sigmin(M(A-\lam I)) = \sigmin(A-\lam I) =
\eta(\lam, A).$ Now defining $ \Delta{A} := - \mu_n M^{-1} U(:,
n)U(:, n)^H,$ we have $\lam \in \Lambda(A+\Delta{A}),\,\, \Delta{A}
\in \S$ and $\|\Delta{A}\|_2 = \eta(\lam, A) = \eta^{\S}(\lam, A).$
Hence the result follows. The proof is similar for the case when $
MA\in\skherm$ and $M^H = -M.$

Finally, consider the case  $\S :=\L$ when $M^H =\pm M.$ Then $MA\in\skherm$ when $ M^H=M$ and $MA\in\herm$ when $M^H=
-M.$ First consider  $MA\in\skherm$ when $ M^H=M.$ Then
for $ \lam \in i\R,$ the set of purely imaginary numbers, we have
$M(A-\lam I) \in \S,$ that is, $M(A-\lam I)$ is skew Hermitian.
Hence the result follows from spectral decomposition of $M(A-\lam
I).$ The proof is similar for the case when $MA\in\herm$ and $M^H=
-M.$ $\blacksquare$

\vone
\noindent{\bf Conclusion.} We have provided a complete solution of the structured mapping problem (Theorem~\ref{gsmp}) for certain classes of structured matrices. More specifically, given a pair of matrices  $X$ and $B$ in $\K^{n\times p}$ and a class $\S \subset \K^{n\times n}$ of structured matrices,  we have provided a complete characterization of structured solutions of the matrix equation $AX =B$ with $A \in \S$ for the case when $\S$ is either a Jordan algebra or a Lie algebra associated with an orthosymmetric scalar product on $\K^n.$ We have determined all optimal solutions in $\S,$ that is, structured solutions which have the smallest norm. We have shown that optimal solution is unique for the Frobenius norm and that there are infinitely many optimal solutions for the spectral norm. We have shown that the results in~\cite{fran6} obtained for a pair of vectors follow as special cases of our general results. Finally, as an application of the structured mapping problem, we have analyzed structured backward errors of approximate invariant subspaces, approximate eigenpairs, and structured pseudospectra of structured matrices.


\begin{thebibliography}{abcd}



\bibitem{Ba:thesis} {\sc B. Adhikari,} {\em Backward perturbation and sensitivity analysis of structured
polynomial eigenvalue problem,}
{PhD thesis, Department of Mathematics, Indian Institute of
Technology Guwahati, India, 2008.}

\bibitem{AlaBKMM} {\sc R. Alam, S. Bora, M. Karow, V. Mehrmann, and J. Moro,} {\em Perturbation theory
for Hamiltonian matrices and the distance to bounded-realness,} { SIAM J. Matrix Anal.
Appl., 32(2011), pp.484 - 514.}


\bibitem{anto}{\sc A. C. Antoulas,} {\em Approximation of Large-Scale Dynamical Systems,}{ SIAM, Philadelphia, 2005.}



%






%

%


\bibitem{D} {\sc C. Davis, W. M. Kahan and H.F. Weinberger,
} {\it Norm-preserving dialations and their applications to
optimal error bounds,} { SIAM J. Numer. Anal., 19(1982), pp.445-469.}








%
%
%


%

\bibitem{graillat:pseudo} {\sc S. Graillat,} {\it A note on structured
pseudospectra,} {J. Comp. Appl. Math., 191(2006) pp. 68-76.}



\bibitem{griv}{\sc S. Grivet-Talocia,}{\em  Passivity enforcement via perturbation of Hamiltonian matrices,}{ IEEE Trans.
Circuits Syst. I. Regul. Pap., 51 (2004), pp. 1755 - 1769.}

\bibitem{horn:book} {\sc R.A. Horn and C. R. Johnson,} {\it Topics in Matrix
Analysis,} {Cambridge University Press, reprinted paperback
version ed., 1995.}

%



\bibitem{karow}{ \sc M. Karow,} {\em $\mu$-values and spectral value sets for linear perturbation classes defined by a scalar
product,} { SIAM. J. Matrix Anal. Appl., 32 (2011), pp. 845 - 865.}






%
%


\bibitem{fran6} {\sc  D. S. Mackey,  N. Mackey and F. Tisseur,} {\it Structured
mapping problems for matrices associated with scalar products,
Part I: Lie and Jordan algebras,} {SIAM J. Matrix Anal. Appl.,
29(2008), pp.1389-1410.}



\bibitem{khatri} {\sc C. Khatri and S. K. Mitra,} {\em  Hermitian and nonnegative definite solutions of linear matrix
equations.,} { SIAM J. Appl. Math., 31 (1976), pp. 579 - 585.}





\bibitem{meinguet} {\sc J. Meinguet,} {\it On the Davis-Kahan-Weinberger solution
of the norm-preserving dilation problem,} {Numer. Math., 49(1986),
pp. 331-341.}



\bibitem{rump}  {\sc S. M. Rump,} {\it  Eigenvalues, pseudospectrum and structured
perturbation,} {  Linear Algebra Appl., 413(2006), pp.567-593.}






%






\bibitem{sun93} {\sc J. Sun,} {\it Backward perturbation analysis of certain
characteristic subspaces,} {Numer. Math., 65 (1993), pp. 357-382.}


\bibitem{tiss03} {\sc F. Tisseur,} { \em A chart of backward errors and condition numbers for singly and doubly structured
eigenvalue problems,}{ SIAM J. Matrix Anal. Appl., 24 (2003), pp. 877 - 897.}

 \bibitem{lntbook} {\sc L. N. Trefethen and M. Embree} { \it  Spectra
and Pseudospectra: the behaviour of nonnormal matrices and
operators,} { Princeton University Press, 2005.}












\end{thebibliography}
\end{document}